\documentclass{TAC}

\usepackage[usenames,dvipsnames]{pstricks}

\usepackage{amscd,amssymb,amsmath,graphicx,verbatim,hyperref}
\usepackage[OT1,T1]{fontenc}

\usepackage[all, cmtip]{xy}
\usepackage{pdfpages}

\usepackage{amsmath,amsfonts,mathabx,tikz,tikz-cd,stmaryrd,url,mathtools, mathrsfs, tensor, pigpen}
\usepackage{hyperref}
\hypersetup{
	colorlinks,
	citecolor=black,
	filecolor=black,
	linkcolor=black,
	urlcolor=black
}

\usetikzlibrary{matrix,arrows,decorations.pathmorphing}

\newcommand{\po}{\ar@{}[dr]|{\text{\pigpenfont R}}}
\newcommand{\pb}{\ar@{}[dr]|{\text{\pigpenfont J}}}
\newcommand{\fpb}{\ar@{}[dr]|{\text{fakepb}}}

\newcommand{\ra}{\rightarrow}
\newcommand{\lra}{\longrightarrow}
\newcommand{\Ra}{\Rightarrow}

\newcommand{\ldual}[1]{\mathord{{\let\nolimits\relax\sideset{^\wedge}{}{#1}}}}
\newcommand{\laction}[2]{\mathord{{\let\nolimits\relax\sideset{^{#1}}{}{#2}}}}
\newcommand{\conj}[2]{\mathord{{\let\nolimits\relax\sideset{^{#1}}{}{#2}}}}

\newcommand{\ox}{\otimes}

\newcommand{\xra}{\xrightarrow}
\newcommand{\xla}{\xleftarrow}

\def\CA{{\mathscr A}}
\def\CB{{\mathscr B}}
\def\CC{{\mathscr C}}

\def\CE{{\mathscr E}}
\def\CF{{\mathscr F}}
\def\CG{{\mathscr G}}

\def\CM{{\mathscr M}}

\def\CP{{\mathscr P}}
\def\CQ{{\mathscr Q}}

\def\CV{{\mathscr V}}

\def\CX{{\mathscr X}}

\makeatletter
\newcommand*\bigcdot{\mathpalette\bigcdot@{.5}}
\newcommand*\bigcdot@[2]{\mathbin{\vcenter{\hbox{\scalebox{#2}{$\m@th#1\bullet$}}}}}
\makeatother

\newcommand{\twocong}[2][0.5]{\ar@{}[#2] \save ?(#1)*{\cong}\restore}
\newcommand{\twoeq}[2][0.5]{\ar@{}[#2] \save ?(#1)*{=}\restore}
\newcommand{\ltwocell}[3][0.5]{\ar@{}[#2] \ar@{=>}?(#1)+/r 0.2cm/;?(#1)+/l 0.2cm/^{#3}}
\newcommand{\rtwocell}[3][0.5]{\ar@{}[#2] \ar@{=>}?(#1)+/l 0.2cm/;?(#1)+/r 0.2cm/^{#3}}
\newcommand{\utwocell}[3][0.5]{\ar@{}[#2] \ar@{=>}?(#1)+/d  0.2cm/;?(#1)+/u 0.2cm/_{#3}}
\newcommand{\dtwocell}[3][0.5]{\ar@{}[#2] \ar@{=>}?(#1)+/u  0.2cm/;?(#1)+/d 0.2cm/^{#3}}
\newcommand{\ultwocell}[3][0.5]{\ar@{}[#2] \ar@{=>}?(#1)+/dr  0.2cm/;?(#1)+/ul 0.2cm/^{#3}}
\newcommand{\urtwocell}[3][0.5]{\ar@{}[#2] \ar@{=>}?(#1)+/dl  0.2cm/;?(#1)+/ur 0.2cm/^{#3}}
\newcommand{\dltwocell}[3][0.5]{\ar@{}[#2] \ar@{=>}?(#1)+/ur  0.2cm/;?(#1)+/dl 0.2cm/^{#3}}
\newcommand{\drtwocell}[3][0.5]{\ar@{}[#2] \ar@{=>}?(#1)+/ul  0.2cm/;?(#1)+/dr 0.2cm/^{#3}}

\begin{document}
\title{The core groupoid can suffice}
\author{Ross Street}
\address{Department of Mathematics, Macquarie University, NSW 2109, Australia}
\eaddress{<ross.street@mq.edu.au>}
\thanks{The author gratefully acknowledges the support of Australian Research Council Discovery Grant DP190102432.}

\date{\small{\today}}
\maketitle
\begin{abstract}
This work results from a study of Nicholas Kuhn's paper entitled 
``Generic representation theory of finite fields in nondescribing characteristic''.
Our goal is to abstract the categorical structure required to obtain an
equivalence between functor categories $[\CF,\CV]$ and $[\CG,\CV]$
where $\CG$ is the core groupoid of the category $\CF$ and $\CV$
is a category of modules over a commutative ring. 
\end{abstract} 
 \tableofcontents
 
 \section*{Introduction}
 
 Let $\mathbb{F}$ be a finite field, let $\mathscr{F}$ be the category of finite dimensional vector spaces 
and linear functions over $\mathbb{F}$, and let $\mathscr{G}$ be the groupoid core of $\mathscr{F}$; that is, the subcategory of $\mathscr{F}$ with all the objects and only the bijective linear functions. 
Let $\mathscr{V}$ be the category of vector 
spaces over the complex numbers (say). In \cite{49} Andr\'e Joyal and I studied a braided monoidal
structure on the functor category $[\mathscr{G},\mathscr{V}]$. 
So I was quite interested and surprised when Nicholas Kuhn told Steve Lack 
and me about his equivalence of categories $[\mathscr{F},\mathscr{V}] \simeq [\mathscr{G},\mathscr{V}]$ 
(see \cite{Kuhn}) when we were working on Dold-Kan-type theorems (see \cite{120} and the related \cite{Helmstutler}). 
The proof of equivalence refers to
a remarkable piece of linear algebra \cite{Kovacs} by our late Australian National University colleague Laci Kov\'acs who built on results in the papers \cite{Faddeev1976, OP1991}. 

Given Kov\'acs' result, the present note is a study of the rest of the argument in order to
abstract the essential categorical features to cover further examples. 

Benjamin Steinberg has recently produced two alternative proofs of the result in \cite{Kovacs}; see
his excellent book \cite{Steinberg2016} and paper \cite{Steinberg2022}. 

Let $\CV$ be the monoidal category of modules over a commutative ring $R$.
Generally, we are interested in categories $\CF$ for which there is a groupoid $\CG$ such that the
functor categories $[\CF,\CV]$ and $[\CG,\CV]$ are equivalent. 
In particular, $\CG$ could be the core groupoid of $\CF$; 
that is, the subcategory with the same objects and with only the invertible morphisms.
Every category $\CF$ gives rise to a $\CV$-category (that is, an $R$-linear category),
denoted $R\CF$, with the same objects and with hom $R$-module $R\CF(A,B)$ free on the
homset $\CF(A,B)$. Indeed, $R\CF$ is the free $\CV$-category on $\CF$ so that the $\CV$-functor category
$[R\CF,\CV]$ is isomorphic to the ordinary functor category $[\CF,\CV]$ with the pointwise $R$-linear structure.
In these terms, we are interested in when $R\CF$ and $R\CG$ are Morita equivalent $\CV$-categories.

In my joint work \cite{120} with Steve Lack on Dold-Kan-type equivalences, we had many examples
of this phenomenon. 
In Section 5 and Appendix A of that work we showed how to reduce our general setting to
the core groupoid case. 
However, the example of Nick Kuhn, where $\CF$ is the category of
finite vector spaces over a fixed finite field $\mathbb{F}$ with all $\mathbb{F}$-linear functions and $\CG$ is the
general linear groupoid over $\mathbb{F}$, does not fit our theory. Yet the ``kernel'' of the equivalence
is of the same type. 
The present work shows that the category theory behind the Kuhn
result also covers Dold-Kan-type results.  

I am grateful to Alexander Campbell for pointing to two gaps in the version of this note that I presented to the Australian Category Seminar on 13 July 2022. I am grateful to Nicholas Kuhn for his help during the writing of this material. Thanks also to Benjamin Steinberg for his valuable feedback on the first arXiv version of this work.

  \section{Review of factorization systems}
  Let $\CF$ be a category equipped with a factorization system $(\CE,\CM)$; see \cite{FreydKelly, KorosTholen}.
  That is, $\CE$ and $\CM$ are sets of morphisms of $\CF$ satisfying
  \begin{itemize}
  \item[FS0.] $mw\in \CM$ if $m\in \CM$ and $w$ is invertible, while $we\in \CE$ if $e\in \CE$ and $w$ is invertible;
  \item[FS1.] if $mh = ke$ with $e\in \CE$ and $m\in \CM$ then there exists a unique $\ell$ with $\ell e = h$ and $m\ell = k$;
  \item[FS2.] every morphism $f$ factors as $f=me$ for some $m\in \CM$ and $e\in\CE$.
  \end{itemize}
  It follows that $\CE\cap\CM$ is the set of invertible morphisms and that $\CE$ and $\CM$ are both closed under composition in $\CF$. 
 We identify $\CE$, $\CM$ and $\CG : = \CE\cap\CM$ with the subcategories of $\CF$ containing all the objects of $\CF$
 but only the morphisms in those respective sets. The following result is a well-known way to express usefully the uniqueness of factorization up to isomorphism.
 \begin{lemma}\label{uniquefactor} Composition in $\CF$ induces an isomorphism
 \begin{eqnarray*}
\int^{C\in \CG}{\CM(C,B)\times \CE(A,C)} \cong \CF(A,B)
\end{eqnarray*}
which is natural in $A\in \CE$ and $B\in \CM$.
 \end{lemma}  
 The factorization system is called {\em proper} when
  \begin{itemize}
  \item[FSP.] every member of $\CE$ is an epimorphism and every member of $\CM$ is a monomorphism.  
    \end{itemize}
   \begin{proposition}\label{propercancellation} In a proper factorization system, if $hf\in\CM$ then $f\in\CM$.
   Dually, if $fk\in\CE$ then $f\in\CE$. 
    \end{proposition} 
     \begin{proof} Assume $hf\in\CM$ and use FS2 to factorize $f= me$ with $e\in \CE$ and $m\in\CM$.
   Then $(hf)1=(hm)e$ so, by FS1, there exists $w$ with $we=1$ and $hfw=hm$. So $e$ is a split monomorphism.
   Using properness, we know $e$ is an epimorphism and hence invertible. So $e\in \CM$, and $f= me\in \CM$.   
    \end{proof}    
    
 \section{Preliminaries in the bicategory of modules}
 Let $\CV$ denote the monoidal category $\mathrm{Ab}^R$ of modules over the commutative ring $R$. 
 
Let $\mathfrak{R}$ denote the bicategory $\CV\text{-}\mathrm{Mod}$. The objects are $\CV$-categories (also called $R$-linear categories and to be thought of as ``$R$-algebras with several objects'' in as much as a category is a ``monoid with several objects'').
The homcategories of $\mathfrak{R}$ are defined to be the $\CV$-functor $\CV$-categories 
$$\mathfrak{R}(\CA,\CB) = [\CB^{\mathrm{op}}\otimes \CA,\CV] \ ; $$
objects of these homs are called {\em modules}.
Composition $$\mathfrak{R}(\CB,\CC)\otimes \mathfrak{R}(\CA,\CB)\xra{\circ}\mathfrak{R}(\CA,\CC)$$
is defined by  
$$(K\circ H)(C,A)=\int^B{H(B,A)\otimes K(C,B)}  \ .$$
The identity module $\CA \to \CA$ is the $\CV$-valued hom functor of $\CA$.

The homcategories of $\mathfrak{R}$ are all abelian $R$-linear categories. 
Consequently, within them, exact sequences have meaning.

For $\CV$-functors $\CA\xra{F}\CC\xla{G}\CB$, we have the module $\CC(G,F) : \CA\to \CB$ with components
$\CC(G,F)(B,A)= \CC(GB,FA)$. In particular, for each $\CV$-functor $F : \CA\to \CB$ we have the module $F_* = \CB(1,F) : \CA\to \CB$. Also, the module $F^*= \CB(F,1) : \CB\to \CA$ provides a right adjoint $F_*\dashv F^*$ for $F_*$ in 
the bicategory $\mathfrak{R}$. 

Here is a result on Morita equivalence for $R$-linear categories; for example, see \cite{50, 88}.

   \begin{proposition}\label{Morita} Let $\CQ \CA$ denote the Cauchy completion of the $\CV$-category $\CA$; that is, the closure of the representables in 
   $[\CA^{\mathrm{op}},\CV]$ under finite direct sums and splittings of idempotents. The following
   conditions on $R$-linear categories $\CA$ and $\CB$ are logically equivalent:
   \begin{itemize}
  \item[(i)] $\CA$ and $\CB$ are equivalent in the bicategory $\mathfrak{R}$;
  \item[(ii)] $[\CA^{\mathrm{op}},\CV] \simeq [\CB^{\mathrm{op}},\CV]$;
  \item[(iii)] $[\CA,\CV] \simeq [\CB,\CV]$;
  \item[(iv)] $\CQ\CA \simeq \CQ\CB$.
  \end{itemize}
    \end{proposition} 
 
We will need Lemma~\ref{retractofadjoint} in our bicategory $\mathfrak{R}$.
A string proof is quite appealing. The result will be used in the proof of Lemma~\ref{idempotentmates}.
 
\begin{lemma}\label{retractofadjoint}
In any bicategory $\mathfrak{M}$, suppose $f\dashv u : A\to X$ is an adjunction with counit $\varepsilon : fu\Ra 1_A$ and unit $\eta : 1_X \Ra uf$. Suppose $\omega : f \Ra f$ is an idempotent 2-morphism on $f$ with splitting provided by $\sigma : f\Ra g$
and $\rho : g\Ra f$: that is, $\rho \sigma = \omega$ and  $\sigma \rho = 1_g$. The mate $\tilde{\omega} : u\Ra u$ of
$\omega : g\Ra g$ under $f\dashv u$ is an idempotent 2-morphism on $g$ satisfying $\varepsilon (\omega u) = \varepsilon (f \tilde{\omega})$. Then, any splitting of $\tilde{\omega}$ delivers a right adjoint $v$ for $g$. Explicitly, if
$\tilde{\rho} \tilde{\sigma} = \tilde{\omega}$ and  $\tilde{\sigma} \tilde{\rho} = 1_v$ then $\tilde{\varepsilon} = \varepsilon (\rho\tilde{\rho})$ and $\tilde{\eta} = (\tilde{\sigma} \sigma)\eta$ provide a counit and unit for $g\dashv v$.
\end{lemma}  
 
 We write $RX$ for the free $R$-module with basis the set $X$. This defines the object assignment of the strong-monoidal left adjoint $R : \mathrm{Set}\to \CV$ to the forgetful functor. 
 We obtain a 2-functor\footnote{The notation of Eilenberg-Kelly \cite{EilKel1966} would be $R_*$ for this 2-functor however we would rather overwork the symbol $R$ than the subscript $*$.}
$R : \mathrm{Cat}\to \CV\text{-} \mathrm{Cat}$ taking ordinary categories to $R$-linear categories
by applying $R$ on the homsets of the categories. The ordinary functor category
\begin{eqnarray*}
[\CF,\CV]
\end{eqnarray*}
with its pointwise $R$-linear enrichment is isomorphic to the $\CV$-enriched $\CV$-functor category
\begin{eqnarray*}
[R\CF,\CV] = \mathfrak{R}(R\CF,R1) \ ,
\end{eqnarray*}
where $1$ is the terminal object of $\mathrm{Cat}$.

Every module $H : \CA\to \CB$ is the ``kernel'' of a functor $\widehat{H} : [\CB,\CV]\to [\CA,\CV]$ defined by
composition with $H$ in $\mathfrak{R}$ thus:
$\widehat{H}(T) = T\circ H : \CA\to R1$. That is,
\begin{eqnarray}\label{hatting}
\widehat{H}(T)A = \int^BH(B,A)\otimes TB \ .
\end{eqnarray}
We obtain a pseudofunctor $\widehat{(-)} : \mathfrak{R}^{\mathrm{op}}\to \CV\text{-}\mathrm{Cat}$.
Moreover, $\widehat{H}$ has a right adjoint $\widetilde{H} : [\CA,\CV]\to [\CB,\CV]$ given by right extension
along $H$ in $\mathfrak{R}$. Explicitly,
\begin{eqnarray}\label{tildeing}
\widetilde{H}(F)B = \int_A[H(B,A), FA] = [\CA,\CV](H(B,-),F) \ . 
\end{eqnarray}

Return now to $\CF$ and its proper factorization system $(\CE,\CM)$ with $\CG=\CE\cap\CM$.
We write $\CM'$, $\CE'$, $\CG'$ for the set of morphisms of $\CF$ not in $\CM$, $\CE$, $\CG$, respectively. 

The $\CV$-functor $j : R\CG\to R\CF$ induced by the inclusion $\CG\hookrightarrow \CF$
gives the right adjoint $\CV$-module $j^*: R\CF \to R\CG$.
So $j^*(A,B) = R\CF(jA,B)$ with $A\in R\CG$ and $B\in R\CF$. 

 \begin{lemma}\label{onM'} For $g\in \CG(C,A)$ and $h\in \CF(B,D)$, there is a commutative square
 \begin{eqnarray}\label{commsq}
 \begin{aligned}
\xymatrix{
R\CM'(A,B) \ar[rr]^-{\subseteq} \ar[d]_-{} && R\CF(A,B) \ar[d]^-{R\CF(g,h)} \\
R\CM'(C,D) \ar[rr]_-{\subseteq} && R\CF(C,D)}
\end{aligned}
\end{eqnarray}
so that $M'(A,B) = R\CM'(A,B)$ defines a submodule $M'$ of $j^*$.
There is a short exact sequence
 \begin{eqnarray}\label{Mses}
\xymatrix{
0 \ar[r]^-{} & M' \ar[rr]^-{\subseteq}  && j^* \ar[rr]^-{q} && M  \ar[r]^-{} & 0 
}
\end{eqnarray}
in $\mathfrak{R}(R\CF,R\CG)$ where $M(A,B)= R\CM(jA,B)$.
Dually, there is a short exact sequence
 \begin{eqnarray}\label{Eses}
\xymatrix{
0 \ar[r]^-{} & E' \ar[rr]^-{\subseteq}  && j_* \ar[rr]^-{\tilde{q}} && E  \ar[r]^-{} & 0 
}
\end{eqnarray}
in $\mathfrak{R}(R\CG,R\CF)$ where $E(B,A)= R\CE(B,jA)$. 
 \end{lemma} 
 \begin{proof}
 Contrapose the implications $hfg\in \CM \ \Ra \ hf = hfgg^{-1}\in\CM \ \Ra \ f\in\CM$ to prove the first sentence.
 The $R$-linear morphism $q_{A,B} : R\CF(A,B)\to R\CM(A,B)$, defined by $q_{A,B}(f) = f$ for $f\in \CM(A,B)$ and $q_{A,B}(f) = 0$ otherwise, has kernel $R\CM'(A,B)$. Then the commutative square \eqref{commsq} 
 induces the module structure on $M$ and renders $q$ a module morphism.
 \end{proof} 
 
 \begin{lemma}\label{ME=1} The family of $R$-module morphisms
 \begin{eqnarray*}
\phi_B : E(B,D) \otimes M(C,B) \to R\CG(C,D) \ ,
\end{eqnarray*}
defined by
$$
\phi_B(e\otimes m) = \begin{cases}
               em & \text{for } em\in \CG \\
0 &\text{otherwise} 
\end{cases}
$$
for $e\in \CE(B,D)$ and $m\in \CM(C,B)$,
is dinatural in $B\in \CF$, natural in $C,D\in \CG$, and induces an invertible morphism $\bar{\phi} : M\circ E\xra{\cong} 1_{R\CG}$ in $\mathfrak{R}$. 
 \end{lemma} 
 \begin{proof}
 Dinaturality requires $\phi_B(ef\otimes m)=\phi_{B'}(e\otimes fm)$ for all $e\in \CE(B',D)$, $m\in \CM(C,B)$ and $f\in \CF(B,B')$.
 This is true since both sides are equal to $efm$ when this is invertible and $0$ otherwise.
 Naturality in $C,D\in \CG$ is obvious. 
 So a natural family $\bar{\phi} : \int^BE(B,D)\otimes M(C,B)\to R\CG(C,D)$ is induced.
 An inverse takes $x\in R\CG(C,D)$ to the equivalence class of $x\otimes 1_C$; it is obviously a right inverse and also
 a left inverse since any $e\otimes m$ representing a non-zero equivalence class has $e\otimes m \sim em\otimes 1_C$
  with $em\in \CG$ by dinaturality of the coend inclusions. 
 \end{proof} 
 
  \begin{corollary}\label{corME=1} The isomorphism $\bar{\phi}$ induces an isomorphism $\widehat{E}\circ \widehat{M}\cong 1_{[\CG,\CV]}$. The mate of this isomorphism, under the adjunction $\widehat{E}\dashv \widetilde{E}$, yields
  a natural transformation $\Theta : \widehat{M}\to \widetilde{E}$. 
  \end{corollary} 
  
  \begin{corollary}\label{corcorME=1} Suppose $\Theta$ too is invertible. Then
  \begin{itemize}
  \item[(i)] $\widehat{E}\dashv \widehat{M}$; 
  \item[(ii)] $\widehat{M}$ is fully faithful; 
  \item[(iii)] $\bar{\phi} : M\circ E\xra{\cong} 1_{R\CG}$ is the counit of an adjunction $M\dashv E$ in $\mathfrak{R}$;
   \item[(iv)] $\widetilde{E}\dashv \widetilde{M}$.
   \end{itemize}   
  \end{corollary} 
\begin{proof} Since $\widehat{E}\dashv \widetilde{E}$ and $\widehat{M}\cong \widetilde{E}$, we obtain (i).   
Since $\widehat{E}\circ \widehat{M}\cong 1_{[\CG,\CV]}$ is the invertible counit of the adjunction in (i), we obtain (ii). 
The pseudofunctor $\widehat{(-)} : \mathfrak{R}^{\mathrm{op}}\to \CV\text{-}\mathrm{Cat}$ is a biequivalence
onto the full subbicategory of $\CV\text{-}\mathrm{Cat}$ consisting of the $\CV$-categories of the form $[\CA,\CV]$ and
left adjoint $\CV$-functors between them; so $\widehat{E}\dashv \widehat{M}\dashv \widetilde{M}$ implies (iii).
Finally, (iv) is obtained by applying the pseudofunctor $\widetilde{(-)} : \mathfrak{R}^{\mathrm{co}}\to \CV\text{-}\mathrm{Cat}$ to $M\dashv E$.   
\end{proof}
    
 There are alternative formulas for $\widehat{M}$ and $\widetilde{E}$ which can be useful. 
 For $A\in \CF$, let $\mathrm{Sub}A$ denote a representative set of the isomorphism classes of $\CM/A$,
 and let $\mathrm{Quo}A$ denote a representative set of the isomorphism classes of $A/\CE$.
 Then
 \begin{eqnarray}\label{altformulas}
\widehat{M}(T)A \ \cong \sum_{(U\xra{m}A)\in \mathrm{Sub}A} TU  & \text{  and  } & \widetilde{E}(T)A \ \cong \prod_{(A\xra{e}W)\in \mathrm{Quo}A} TW \ . 
\end{eqnarray}

In cases where $\mathrm{Sub}A$ and $\mathrm{Quo}A$ are finite, both the sum and product in \eqref{altformulas}
are direct sums and the components of $\Theta : \widehat{M}\to \widetilde{E}$ transport across the isomorphisms 
to matrices whose non-zero entries 
have the form $T(em) : TU\to TW$ for $em\in \CG$. Suppose that elements of each $\mathrm{Sub}A$ can be listed 
$m_0, \dots, m_s$ and those of each $\mathrm{Quo}A$ can be listed $e_0, \dots, e_s$, in such a way that $e_im_j\in \CG$ implies $1\le i \le j \le s$ and that $e_im_i=1$ for $0\le i \le s$. Then these matrices are square and triangular 
with identities down the main diagonal, 
and so (because we have additive inverses in the hom $R$-modules) are invertible.  

 \begin{proposition}\label{locpo} Suppose $\CF$ is the underlying category of a category $\CP$ enriched in the
 cartesian monoidal category of partially ordered sets; so $\CP$ is a locally ordered 2-category. 
 Suppose each $m\in \CM$ has a right adjoint $m^*\in \CE$ in $\CP$ with identity unit.
 Suppose each object of $\CF$ has only finitely many automorphisms. 
 Then a reflexive, antisymmetric, transitive relation $m\trianglelefteq n$ is defined on $\mathrm{Sub}A$
 by: there exists $x\in \CG$ with $mx \le n$. 
  \end{proposition} 
  \begin{proof} This is a variant of the proof of Proposition 2.9 in \cite{120}.
  
  Reflexivity is clear since each $1\in \CG$. 
  
  Suppose $m\trianglelefteq n \trianglelefteq \ell$.
  Then $mx \le n$ and $ny \le \ell$ for some $x,y\in \CG$.
  So $mxy \le ny\le \ell$ yields $m\trianglelefteq \ell$ and so transitivity. 
  
  Suppose $m\trianglelefteq n \trianglelefteq m$. Then $mx \le n$ and $ny \le m$ for some $x,y\in \CG$. 
  So $mxy \le ny\le m$. Since $m$ is fully faithful in $\CP$, $xy\le 1$ and we have the chain
  $$\dots \le (xy)^r\le (xy)^{r-1}\le \dots \le (xy) \le 1$$
  which cannot be infinite since the domain of $m$ has only finitely many automorphisms.
  Since $xy$ is invertible, this means $(xy)^r=1$ for some $r \ge 0$.
  Then $1 = (xy)^r \le \dots \le (xy) \le 1$ yielding $xy = 1$. So $m = mxy \le ny\le m$ implies $m = ny$.
  Thus $m$ and $n$ are in the same isomorphism class of $\CM/A$; since they are both in $\mathrm{Sub}A$,
  $m=n$ and we have antisymmetry.   
  \end{proof}    

 \begin{corollary}\label{corlocpo} In the situation of Proposition~\ref{locpo}, suppose $\mathrm{Sub}A$ is finite
 and every $(A\xra{e}W)\in \CE$ is isomorphic to $m^*$ for some $m\in \mathrm{Sub}A$.
 Then there is a listing of $m_0, \dots, m_s$ of the elements of $\mathrm{Sub}A$ such that, together with
 the listing $m_0^*, \dots, m_s^*$ of $\mathrm{Quo}A$, gives an invertible matrix representing 
 $\Theta : \widehat{M}\to \widetilde{E}$.       
  \end{corollary} 
  \begin{proof} Take the linear order $m_0, \dots, m_s$ to contain the order $\trianglelefteq$ of the Proposition.
  Then $x = m_i^*m_j \in \CG$ implies $m_i x \le m_j$ and so $m_i\trianglelefteq m_j$ and then $i\le j$.    
  \end{proof}
    
  \begin{example}\label{exlocpo} Let $\CF = \Delta_{\bot}$ be the category of finite non-empty ordinals 
  $\mathbf{a} = \{0,1,\dots, a-1\}$ and order-preserving functions $\xi : \mathbf{a}\to \mathbf{b}$ which 
  also preserve first element. Let $(\CE,\CM)$ be the surjective-injective factorization system. Let $\CP$
  be the 2-category obtained by ordering the homs of $\CF$ with the pointwise order. Each morphism
  $\xi : \mathbf{a}\to \mathbf{b}$ of $\CF$ has a right adjoint $\xi^*$ as a functor since $\xi$ preserves
  the initial object. If $\xi\in\CM$ then $\xi^*\in \CE$ with $1= \xi^*\xi$, so Corollary~\ref{corlocpo} applies.   
  As the only invertible morphisms in $\CF$ are identities, the groupoid $\CG$ is discrete with
  countably many objects; so $[\CG,\CV]$ is the product $\CV\times \CV\times \dots $ of countably
  many copies of $\CV$, that is, the category of objects of $\CV$ graded by the positive integers.        
  \end{example} 

 \begin{example}\label{2exlocpo} 
 Let $\CA$ be any category in which every morphism is a monomorphism, pullbacks exist,
 and each slice category $\CA/A$ has finitely many isomorphism classes.
 Take $\CF$ to be the category of spans in $\CA$; that is, the objects are those of $\CA$ and
 the morphisms $f : A\to B$ are isomorphism classes $[f_1,U,f_2]$ of spans $A\xla{f_1}U\xra{f_2}B$ in $\CA$,
 where composition uses pullback. Take $\CM$ to consist of those spans $f$ with $f_1$ invertible and
 $\CE$ to consist of those spans $f$ with $f_2$ invertible. Let $\CP$ be the 2-category obtained by
 declaring $f = [f_1,U,f_2]\le [g_1,V,g_2]=g : A\to B$ when there exists $h : U\to V$ in $\CA$ with $g_1h=f_1$
 and $g_2h=f_2$. Corollary~\ref{corlocpo} applies. The core groupoid $\CG$ of $\CF$ is isomorphic to the
 core groupoid of $\CA$. An example of this is $\CA= \mathrm{FI}$, the category of finite sets and injective
 functions. Then $\CG = \mathfrak{S}$, the groupoid of finite sets and bijections, sometimes called
 the symmetric groupoid. So $[\CG,\CV]$ is the category of $R$-linear Joyal species \cite{AnalFunct};
 it is equivalent to the product of the categories of $R$-linear representations of the symmetric groups
  $\mathfrak{S}_n$, $n\ge 0$.
  \end{example} 
  
  \begin{remark}\label{rmklocpo} Under the assumptions in Proposition~\ref{locpo} and Corollary~\ref{corlocpo}, 
  the functor $\widehat{M} :[\CG,\CV]\to [\CF,\CV]$ is actually an equivalence; so $M\dashv E$ is
  an adjoint equivalence $R\CF \simeq R\CG$ in $\mathfrak{R}$. 
  A proof is in Section 6 of \cite{120} so we shall not repeat it here. 
  Therefore already we have numerous examples
  where the core suffices for representations of $\CF$ in $\CV$.
  In the next two sections we shall see how every finite field provides an example.
  \end{remark} 

\section{Main theorem}\label{MainThmSection}
We continue with a proper factorization system $(\CE,\CM)$ on a category $\CF$.

The goal is to find other conditions under which $M : R\CF\to R\CG$ in $\mathfrak{R}$ is an equivalence.
If this is to be the case, $M$ will need to have a left adjoint so we would like to apply Lemma~\ref{retractofadjoint}
to deduce the adjoint from the adjunction $j_*\dashv j^*$. This means we would like $M$ to be a retract of $j^*$.
A natural splitting of the short exact sequence \eqref{Mses} would suffice. 

A splitting $\rho : M \to j^*$ of the epimorphism $q : j^* \to M$ has components $\rho_{A,B} : M(A,B)\to R\CF(A,B)$
which are natural in $A\in \CG$ and $B\in \CF$. In particular they are natural in $B\in R\CM$ and so, by Yoneda,
are determined by the value, say $p_A\in R\CF(A,A)$, at $1_A\in R\CM(A,A)$. Then $\rho_{A,B}(m)=mp_A$.
This motivates the following assumption.  
\bigskip

\noindent \textbf{Idempotent Axiom.} For each object $A\in \CF$, there is a morphism $p_A : A \to A$ in $R\CF$ such that 
\begin{itemize}
  \item[p0.] $p_Ap_A=p_A;$
  \item[p1.] if $f\in \CF(A,B)$ and $f \in \CM'$ then $fp_A = 0$;
  \item[p2.] if $f\in \CF(A,B)$ and $f \in \CE'$ then $p_Bf = 0$;
   \item[p3.] if $g\in \CG(A,B)$ then $gp_A = p_Bg$;
    \item[p4.] $p'_A= 1_A- p_A \in R(\CE'(A,A)\cap \CM'(A,A))$.
  \end{itemize}  

\begin{example}
Every groupoid $\CG$ provides a trivial example with $\CF = \CE = \CM = \CG$ and $p_A = 1_A$.
In the general case, we can think of $p'_A=1_A-p_A$ as an obstruction to $\CF$'s being a groupoid
although we will see in Theorem~\ref{mainthm} a sense in which $\CF$ is not too far from $\CG$
at the $R$-linear level.  
\end{example}

Using p1, we have the identity linear function as the left side of a commutative square \eqref{sigma}.
By Proposition~\ref{propercancellation}, if $u\in \CM'(A,A)$ and $f\in \CF(A,B)$ then $fu\in \CM'$; so, 
by p4, $fp'_A\in R\CM'(A,B)$ for all $f\in \CF(A,B)$. 
This gives a splitting $\sigma_{A,B}$ of the idempotent $R\CF(p'_A,1_B)$.
 \begin{eqnarray}\label{sigma}
 \begin{aligned}
\xymatrix{
R\CM'(A,B) \ar[rr]^-{\subseteq} \ar[d]_-{1} && R\CF(A,B) \ar[d]^-{R\CF(p'_A,1_B)} \ar[lld]_-{\sigma_{A,B}}\\
R\CM'(A,B) \ar[rr]_-{\subseteq} && R\CF(A,B)}
\end{aligned}
\end{eqnarray}
Using p3, we see that the linear functions $R\CF(p'_A,1_B)$ are the components of an idempotent module morphism
$\pi' : j^*\ra j^*$ which is split by the module morphism $\sigma : j^* \ra M'$ and the inclusion $M'\hookrightarrow j^*$.
It follows that the short exact sequence \eqref{Mses} splits so that, if we put  $\pi = 1-\pi' = R\CF(p_A,1_B)$, we have:
 \begin{lemma}\label{rhosplitting} 
 The idempotent module morphism $\pi = R\CF(p_A,1_B) : j^*\to j^*$ splits as 
 \begin{eqnarray}\label{splitpi}
 \begin{aligned}
\xymatrix{
j^* \ar[rr]^-{\pi} \ar[rd]_-{q} && j^* \ar[rd]^-{q}  \\
& M \ar[ru]_-{\rho} \ar[rr]_-{1_M} && M}
\end{aligned}
\end{eqnarray}
where, for $m\in\CM(A,B)$, $\rho_{A,B}(m) = mp_A$ and for 
$g \in \CG(A',A)$ and $f\in \CF(B,B')$, 
\begin{eqnarray*}
M(g,f)m = fmgp_{A'} \ .
\end{eqnarray*}
 \end{lemma}
Dually, we also have:
\begin{lemma}\label{tilderhosplitting} 
 The idempotent module morphism $\varpi = R\CF(1_A,p_B): j_*\to j_*$ splits as 
 \begin{eqnarray}\label{splitvarpi}
 \begin{aligned}
\xymatrix{
j_* \ar[rr]^-{\varpi} \ar[rd]_-{\tilde{q}} && j_* \ar[rd]^-{\tilde{q}} \\
& E \ar[ru]_-{\tilde{\rho}}  \ar[rr]_-{1_E} && E}
\end{aligned}
\end{eqnarray}
where $E(A,B) = R\CE(A,B)$, $\tilde{\rho}_{A,B}(e) = p_Be$, and
\begin{eqnarray*}
E(f,g)e = p_{B'}gef \ ,
\end{eqnarray*}
 for $e\in\CE(A,B)$, $f\in \CF(A',A)$, $g\in \CE((B,B')$.
 \end{lemma}
 
 \begin{lemma}\label{idempotentmates} 
 The idempotent module morphisms $\pi : j^*\to j^*$ and $\varpi : j_*\to j_*$ are mates under the adjunction $j_*\dashv j^*$ in $\mathfrak{R}$. Consequently, $E\dashv M$ in $\mathfrak{R}$.
 \end{lemma}
 \begin{proof} 
 By the Yoneda Lemma or using p3 directly, the square
 \begin{equation*}
\xymatrix{
R\CG(C,D) \ar[rr]^-{j} \ar[d]_-{j} && R\CF(jC,jD) \ar[d]^-{R\CF(p_C,1_{jD})} \\
R\CF(jC,jD) \ar[rr]_-{R\CF(1_{jC},p_D)} && R\CF(jC,jD)}
\end{equation*}
commutes and the arrows marked $j$ are the components of the unit of $j_*\dashv j^*$.
This proves the first sentence. The second sentence uses Lemma~\ref{retractofadjoint}. 
 \end{proof}

Since $\widehat{(-)} : \mathfrak{R}^{\mathrm{op}}\to \CV\text{-}\mathrm{Cat}$ (see \eqref{hatting}) is a pseudofunctor, we have:

\begin{corollary}\label{hatadjunction} 
 The functor $\widehat{M} : [\CG,\CV]\to [\CF,\CV]$  
 has right adjoint $\widehat{E}$.
 \end{corollary}
 
 Since $\widetilde{(-)} : \mathfrak{R}^{\mathrm{co}}\to \CV\text{-}\mathrm{Cat}$ (see \eqref{tildeing}), we have:
 
 \begin{corollary}\label{tildeadjunction} 
 The functor $\widetilde{M} : [\CF,\CV]\to [\CG,\CV]$
 has right adjoint $\widetilde{E}$. 
 \end{corollary}
 
 \begin{corollary}\label{EsareMs} 
 From $\widehat{M}\dashv \widetilde{M}$ and $\widehat{E}\dashv \widetilde{E}$, it follows that $\widehat{E}\cong \widetilde{M}$.     
 \end{corollary}
 
 \begin{lemma}\label{invertibleunit} 
 The unit $1_{R\CG}\Ra M\circ E$ of the adjunction $E\dashv M$ in $\mathfrak{R}$ is invertible. So the unit 
 $1_{[\CG,\CV]}\Ra \widehat{E}\circ \widehat{M}$ of $\widehat{M}\dashv \widehat{E}$ is invertible and $\widehat{M}$ is fully faithful.
 \end{lemma}
 \begin{proof} Note that $(j^*\circ j_*)(C,D) \cong R\CF(jC,jD)$ and, from Lemma~\ref{retractofadjoint}, the unit of the adjunction $E\dashv M$
 has component at $(C,D)$ equal to the composite
 \begin{eqnarray*}
R\CG(C,D)\xra{j} R\CF(jC,jD)\cong \int^{B\in \CF}R\CF(B,jD)\otimes R\CF(jC,B)\xra{\int^{B\in \CF}q\otimes \tilde{q}} (M\circ E)(C,D)
\end{eqnarray*}
which is inverse to the component of the isomorphism $\bar{\phi}$ of Lemma~\ref{ME=1}.

Alternatively, we can composing the splittings of Lemmas~\ref{rhosplitting} and \ref{tilderhosplitting} yielding the splitting
  \begin{eqnarray*}
 \begin{aligned}
\xymatrix{
R\CF(jC,jD) \ar[rr]^-{\pi\circ \varpi = R\CF(p_C,p_D)} \ar[rd]_-{q\circ \tilde{q}} && R\CF(jC,jD) \ar[rd]^-{q\circ \tilde{q}}  \\
& (M\circ E)(C,D) \ar[ru]_-{\rho\circ\tilde{\rho}} \ar[rr]_-{1} && (M\circ E)(C,D) \ .}
\end{aligned}
\end{eqnarray*}
However, we also have the factorization
 \begin{eqnarray*}
 \begin{aligned}
\xymatrix{ 
R\CF(jC,jD) \ar[rd]_-{\tilde{q}} \ar[rr]^-{R\CF(1_C,p_D)} && R\CF(jC,jD) \ar[rd]_-{q} \ar[rr]^-{R\CF(p_C,1_D)} && R\CF(jC,jD) \\
& R\CE(jC,D) \ar[rd]_-{q\restriction} \ar[ru]_-{\subseteq} && R\CM(C,jD) \ar[ru]_-{\subseteq} & \\
&& R\CG(C,D) \ar[ru]_-{\subseteq} &&
}
\end{aligned}
\end{eqnarray*}
also giving a splitting of the idempotent $\pi\circ \varpi$. The induced isomorphism $R\CG(C,D)\cong (M\circ E)(C,D)$ 
is the component of the unit for $E\dashv M$ at $(C,D)$.  
\end{proof}

\begin{lemma}\label{retractcounit} 
 If each object $A$ of $\CF$ has only finitely many $\CM$-subobjects (that is, each slice category $\CM/A$ has only
 finitely many isomorphism classes) then the counit of the adjunction $E\dashv M$ in $\mathfrak{R}$ is a split epimorphism.
 \end{lemma}
 \begin{proof} 
 First a non-proof! Since $R : \mathrm{Set}\to \CV$ preserves colimits and tensor products, Lemma~\ref{uniquefactor} yields an
 isomorphism 
 \begin{eqnarray*}
(E\circ M)(A,B) = \int^{C\in \CG}{R\CM(C,B)\otimes R\CE(A,C)} \cong R\CF(A,B)
\end{eqnarray*}
which in natural in $A\in \CE$ and $B\in \CM$, not as objects of $\CF$.
Moreover, as we see from Lemma~\ref{retractofadjoint}, this isomorphism is not the counit $E\circ M\Ra 1_{R\CF}$ of $E\dashv M$. 

Now for the proper proof. The counit, which is natural, is the composite
\begin{eqnarray*}
(E\circ M)(A,B) = \int^{C\in \CG}{R\CM(C,B)\otimes R\CE(A,C)} \ra (j_*\circ j^*)(A,B) \ra R\CF(A,B)
\end{eqnarray*}
which takes the equivalence class of $m\otimes e \in R\CM(C,B)\otimes R\CE(A,C)$ to $mp_Ce \in R\CF(A,B)$.
We claim this family has a natural right inverse. In particular, for $A=B$, we must see that the identity $1_A$ of $A$
is in the image. Take a finite family $(C_i\xra{m_i}A)_{i=0}^k$ of morphisms in $\CM$ representing all isomorphism
classes in the ordered set $\CM/A$ and with the property that $m_i = m_jn$ for some $C_i\xra{n}C_j$ implies $i\le j$;
we can suppose $C_k=A$ and $m_k=1_A$. Let
\begin{eqnarray*}
\phi_j : \bigoplus_{i\le j}{R\CM(C_i,A)\otimes R\CE(A,C_i)} \lra R\CF(A,A)
\end{eqnarray*}
denote the function defined by $\phi_j(m\otimes e) = mp_{C_i} e$ for $m\in R\CM(C_i,A)$ and $e\in R\CE(A,C_i)$.
So the image of $\phi_i$ is contained in the image of $\phi_j$ for $i\le j$.
We will show that every $f\in \CF(A,A)$ is in the image of the function $\phi_k$. 
The proof is by induction on $j$, where $f=m_je$ for some $e\in \CE$ and unique $j\le k$. 
As $C_0$ has no proper $\CM$-subobjects,
$\CM'(C_0,C_0) = \emptyset$. By Idempotent Axiom p4, we have $p'_{C_0} = 0$. So $p_{C_0} = 1_{C_0}$ and 
$f =m_0e = m_0p_{C_0}e = \phi_0(m_0\otimes e)$ is in the image of $\phi_0$. Now take $0< j\le k$ and assume the result for $j-1$. 
Then $f=m_je = m_jp_{C_j}e + m_jp'_{C_j}e = \phi_j(m_j\otimes e) + m_jp'_{C_j}e$.
By Idempotent Axiom p4, $p'_{C_j}$ is an $R$-linear combination of terms of the form $ne'$ with $n\in \CM$ and
$e'\in \CE$ non-invertible. So terms of $m_jp'_{C_j}e$ are of the form $m_ie''$ with $i< j$ and $e''\in \CE$.
By induction, all these terms are in the image of $\phi_{j-1}$. So $f$ is in the image of $\phi_j$. 

Hence, $1_A\in \CF(A,A)$ is in the image of the $A, A$ component of the counit. Suppose $t\in  \int^{C\in \CG}{R\CM(C,A)\otimes R\CE(A,C)}$ maps to $1_A\in \CF(A,A)$. By the Yoneda Lemma, there exists a unique family of morphisms
\begin{eqnarray*}
R\CF(A,B) \lra \int^{C\in \CG}{R\CM(C,B)\otimes R\CE(A,C)}
\end{eqnarray*}
 which is natural in $B\in \CF$ taking $1_A\in \CF(A,A)$ to $t$.
 Again by Yoneda, this gives right inverses to the components of the counit. 
 Dually, there is such a family natural in $B$. By Yoneda yet again, these families agree
 since they have the same value at the identity.
 So the counit of $E\dashv M$ has a natural right inverse.       
 \end{proof}
 
 \begin{theorem}\label{mainthm} 
If each object $A$ of $\CF$ has only finitely many $\CM$-subobjects then the adjunction $E\dashv M$ is an equivalence $R\CF \simeq R\CG$ in $\mathfrak{R}$.
 \end{theorem}
 \begin{proof} 
Using Lemma~\ref{invertibleunit} and Corollary~\ref{EsareMs}, we know that the unit $1_{R\CG}\Ra \widetilde{M}\circ \widehat{M}$
is invertible. So $\widehat{M} : [R\CG,\CV]\to [R\CF,\CV]$ is fully faithful and so is its composite $M' : R\CG^{\mathrm{op}} \to [R\CF,\CV]$, $C\mapsto M(C,-)$, with the Yoneda embedding. Moreover, Lemma~\ref{retractcounit} implies $M'$ is strongly generating. 
By Theorem 2 of \cite{29}, $M'$ is also dense.
So the counit of $E\dashv M$ is invertible. 
 \end{proof}

In the terminology of ring theory in the spirit of \cite{50}, 
Theorem~\ref{mainthm} implies that $R\CF$ and $R\CG$ are Morita equivalent
several-object $R$-algebras. In the terminology of enriched category theory, it implies that $R\CF$ and $R\CG$ are Cauchy equivalent $\CV$-categories. 

\begin{corollary}\label{maincorollary} 
 The functor $\widehat{M} : [\CG,\CV]\to [\CF,\CV]$   
 is an equivalence with inverse equivalence $\widehat{E}$.
 \end{corollary}

 \begin{corollary}\label{maincorollarydual} 
 The functor $\widetilde{M} : [\CF,\CV]\to [\CG,\CV]$
 is an equivalence with inverse equivalence $\widetilde{E}$.
 \end{corollary}
 
  \begin{corollary}\label{mainfunctorassplitting} 
  The equivalence $\widetilde{M}$ is a retract of the restriction functor 
  $$[j,1] : [R\CF, \CV] \to [R\CG, \CV] \ .$$
 \end{corollary}
 \begin{proof} 
 Apply the contravariant functor $[\CF,\CV](-, F)$ to diagram \eqref{splitpi} of Lemma~\ref{rhosplitting} to obtain
 the idempotent splitting
  \begin{eqnarray*}
 \begin{aligned}
\xymatrix{
[\CF,\CV](R\CF(j,-), F) \ar[rr]^-{[\CF,\CV](\pi, F)} \ar[rd]_-{[\CF,\CV](\rho, F)} && [\CF,\CV](R\CF(j,-), F)  \\
& [\CF,\CV](M(1_{R\CG},-), F) \ar[ru]_-{[\CF,\CV](q, F)}  & \ .}
\end{aligned}
\end{eqnarray*} 
The top side of the triangle transports across the Yoneda isomorphism to an idempotent on $Fj$ while the bottom vertex is
isomorphic to $\widetilde{M}$.  
 \end{proof}
 \begin{corollary}\label{generalcodomain} 
  For any $R$-linear category $\CX$ in which has finite direct sums and splitting of idempotents, there is an equivalence 
  $$[\CF, \CX] \simeq [\CG, \CX] \ .$$
 \end{corollary} 
 
 \section{Stiffness and Kov\'acs idempotents}\label{SaKi}
 
 \begin{definition}
In any category $\CF$, a morphism $f : A\to B$ is called {\em stiff} when the only endomorphisms of $A$ through which $f$
factors are automorphisms. In other words, $f = (A\xra{u}A\xra{v}B)$ implies $u$ invertible. A morphism is {\em costiff} 
when it is stiff in the opposite category. A category is called {\em stiff} when the costiff and stiff morphisms are the $\CE$
and $\CM$ of a proper factorization system.
    \end{definition}
The category of finite sets and the category of finite-dimensional vector spaces over a field are stiff categories: in both costiff means surjective and stiff means injective.

Clearly any stiff endomorphism is an automorphism so a consequence of stiffness is what might be called ``the pigeon-hole principle'':
 \begin{proposition}\label{PhP}
For all objects $A$ in a stiff category $\CF$, the inclusions $\CM(A,A)\subseteq \CG(A,A)$ and $\CE(A,A)\subseteq\CG(A,A)$ are equalities. Equivalently, if $\CG(A,B)\neq \emptyset$ then $\CG(A,B)=\CE(A,B)=\CM(A,B)$. 
  \end{proposition}  
  
 Recall that the set $\CG'$ of morphisms is the complement of $\CG$ in $\CF$. 
 For each $C\in \CF$, Proposition~\ref{PhP} implies that $R\CG'(C,C)$ is a two-sided ideal in $R\CF(C,C)$. 
 In particular $R\CG'(C,C)$ is an $R$-algebra, possibly without an identity element.
  In the case where $\CF = \mathrm{vect}_{\mathbb{F}}$ is the category of finite vector spaces over a finite field $\mathbb{F}$ and the factorization
  is surjective-injective $\mathbb{F}$-linear functions, Laci Kov\'acs \cite{Kovacs} showed that there is an identity element making $R\CG'(C,C)$
  a unital algebra provided the characteristic of $\mathbb{F}$ is invertible in $R$;
  also see Subsection 3.1 of Kuhn's paper \cite{Kuhn} for helpful material on this.
  Steinberg has two alternative methods of producing these identity elements; see \cite{Steinberg2016, Steinberg2022}. 
  Therefore we make the following definition in the general situation.
  \begin{definition}
  A {\em Kov\'acs idempotent} is an identity element $\ell_C$ making $R\CG'(C,C)$ an $R$-algebra under composition.  
    \end{definition}
    
    \begin{example}\label{ExDeltabot3}
  Return to Example~\ref{exlocpo} where $\CF = \Delta_{\bot}$ and take $C = \mathbf{3} = \{0,1,2\}$.
  The monoid $R\CF(C,C)= \Delta_{\bot}(\mathbf{3},\mathbf{3})$ has six elements, namely, the functions
  \begin{eqnarray*}
  \iota = 
\begin{bmatrix}
0 & 1 & 2 \\
0 & 1 & 2
\end{bmatrix} \qquad 
 \sigma_0 = 
 \begin{bmatrix}
0 & 1 & 2 \\
0 & 0 & 2
\end{bmatrix}\qquad 
 \sigma_1 = 
 \begin{bmatrix}
0 & 1 & 2 \\
0 & 1 & 1
\end{bmatrix}
\end{eqnarray*}
\begin{eqnarray*}
  \tau = 
\begin{bmatrix}
0 & 1 & 2 \\
0 & 2 & 2
\end{bmatrix} \qquad 
 \upsilon = 
 \begin{bmatrix}
0 & 1 & 2 \\
0 & 0 & 1
\end{bmatrix}\qquad 
 \zeta = 
 \begin{bmatrix}
0 & 1 & 2 \\
0 & 0 & 0
\end{bmatrix}
\end{eqnarray*}  
    A presentation of this monoid is provided by the generators $\sigma_0$, $\sigma_1$ and $\tau$
    subject to the relations that the generators are idempotents and
    \begin{eqnarray*}
\sigma_0\sigma_1\sigma_0 = \sigma_0\sigma_1 = \sigma_1\sigma_0\sigma_1 \phantom{aaaaaaa}\\
\tau \sigma_0 = \sigma_0 \ , \ \tau \sigma_1 = \tau = \sigma_0\tau \ , \ \sigma_1\tau = \sigma_1 \ .
\end{eqnarray*}
The only invertible element is the identity $\iota$. So $\CG'(C,C)$ has the five other elements.
The Kov\'acs idempotent in $R\CG'(C,C)$ is $\ell_{\mathbf{3}} = \sigma_0 + \sigma_1 -\upsilon$.    
    \end{example}
\bigskip
  
  Let us clarify some terminology for use in upcoming Remarks~\ref{RmkDeltabot3} and \ref{RmkKovacsId}. 
For an object $C$ in a category $\CF$, we write $\CF(C)$ for the full subcategory of $\CF$
with the one object $C$. With the viewpoint that monoids are one object categories, $\CF(C)$ is identified with the monoid $\CF(C,C)$.

 \begin{remark}\label{RmkDeltabot3} 
Perhaps this is a good point to add a little to Example~\ref{ExDeltabot3} because
it does provide an example where a groupoid suffices.
The core groupoid of $\Delta_{\bot}(\mathbf{3})$ is a discrete category with one object.
Clearly $R\Delta_{\bot}(\mathbf{3})$ is not Morita equivalent to $R1$.
However, in $R\Delta_{\bot}(\mathbf{3})$, we have a complete list of
orthogonal idempotents (see Appendix~B of \cite{120}):
\begin{eqnarray*}
e_0 = \sigma_0\sigma_1 \ , \  e_1 = (1-\sigma_0)\sigma_1 \ , \  e_2 = (1-\sigma_1)\sigma_0 \ , \  e_3 = (1-\sigma_1)(1-\sigma_0) \ ,
\end{eqnarray*}
so that, in the Cauchy completion $\CQ R\Delta_{\bot}(\mathbf{3})$ (see Proposition~\ref{Morita}), we have 
\begin{eqnarray*}
\mathbf{3} \cong E_0 \oplus E_1 \oplus E_2 \oplus E_3  
\end{eqnarray*}
 where $E_i$ is the object obtained in splitting the idempotent $e_i$. 
 If we put $\gamma = \sigma_1\sigma_0 - \sigma_0\sigma_1$ and $\delta = \tau -\sigma_1$, using the relations in the presentation, we obtain the equations $\gamma\delta = e_1$ and $\delta\gamma = e_2$ yielding $E_1 \cong E_2$. It follows that we have an isomorphism of the form
 \begin{eqnarray*}
\mathbf{3} \cong A \oplus B \oplus B \oplus C \ .  
\end{eqnarray*}
Transporting the endomorphisms $\sigma_0$, $\sigma_1$, $\tau$ on $\mathbf{3}$ across this isomorphism
yields the matrices 
  \begin{eqnarray*}
\begin{bmatrix}
1 & 0 & 0 & 0 \\
0 & 0 & 1 & 0 \\
0 & 0 & 1 & 0 \\
0 & 0 & 0 & 0
\end{bmatrix} \ , \qquad 
 \begin{bmatrix}
1 & 0 & 0 & 0 \\
0 & 1 & 0 & 0 \\
0 & 0 & 0 & 0 \\
0 & 0 & 0 & 0
\end{bmatrix}  \ ,\qquad 
 \begin{bmatrix}
1 & 0 & 0 & 0 \\
0 & 1 & 0 & 0 \\
0 & 1 & 0 & 0 \\
0 & 0 & 0 & 0
\end{bmatrix}  \ .
\end{eqnarray*}
Any endomorphism $f$ of $A \oplus B \oplus B \oplus C$ which commutes with these three matrices
must be of the form $f = f_0 \oplus f_1\oplus f_1 \oplus f_2$.
It follows that $R\Delta_{\bot}(\mathbf{3})$ is $\CV$-Morita equivalent to $R3$ where $3$
is the discrete category with three objects. In other words, we have an equivalence of $R$-linear categories 
$$[\Delta_{\bot}(\mathbf{3}), \CV]\cong \CV\times \CV \times \CV \ .$$
 \end{remark}
\bigskip
    
    Let us clarify some more terminology for use in the upcoming Remark~\ref{RmkKovacsId}. A morphism $f : A\to B$ in a category is said to {\em split} (or to be {\em von Neumann regular}) when there exists a morphism $g : B\to A$ such that $fgf=f$. A monomorphism splits if and only if it has a left inverse; that is, it is a coretraction. An epimorphism splits if and only if it has a right inverse; that is, it is a retraction. If a morphism $f$ factors as $f=jr$ with $r$ a retraction and $j$ a coretraction then $f$ splits. A category is {\em von Neumann regular} when every morphism is split;
    this terminology applies to monoids as one object categories.  
    
    Every morphism in the category of sets whose domain is not empty splits. The category of vector spaces over a division ring is von Neumann regular as are the categories $\CF$ of Examples~\ref{exlocpo} and \ref{2exlocpo}. 
    
 \begin{remark}\label{RmkKovacsId} 
 It can happen that $R\CF(C)$ is semisimple. In particular, this is true when $Y = \CF(C)$
 is any finite monoid of Lie type and $R$ is a field of characteristic 0; see Corollary 2.9 of \cite{OP1991}. Section 5.4 of Steinberg \cite{Steinberg2016} provides necessary and sufficient conditions on
 a finite monoid $Y$ and field $R$ in order for $RY$ to be semisimple: it is necessary that $Y$ be regular and the characteristic of $R$ should not divide the order of the group of invertible elements
 of the monoid $pYp$ for any idempotent $p\in Y$.  
 With semisimplicity of $RY$, the inclusion of every two-sided ideal $J$ of $RY$ into $RY$ splits
 as a left module morphism and splits as a right module morphism. The value of the left module
 splitting at the identity of $RY$ gives a right identity for $J$ and the value of the right module
 splitting at the identity of $RY$ gives a left identity for $J$. So $J$ becomes an $R$-algebra with
 identity. This is a source of Kov\'acs idempotents for categories $\CF$ whose endomorphism
 monoids $\CF(C)$ are finite and have $R\CF(C)$ semisimple.   
 \end{remark}
 
  \begin{example}\label{KovacsEx}
  As already mentioned, if $\mathbb{F}$ is a finite field, $\CF = \mathrm{vect}_{\mathbb{F}}$, and $R$ is a ring in which the characteristic of $\mathbb{F}$ is invertible then each object $C\in \CF$ has a Kov\'acs idempotent $\ell_C$; see \cite{Kovacs}.  
    \end{example}
    
     \begin{example}\label{OPEx}
Let $\mathbb{F}$ be a finite field. Let $\CF$ be the category whose objects $(V,\alpha)$ are
finite $\mathbb{F}$-vector spaces $V$ equipped with an $\mathbb{F}$-linear isomorphism $\alpha : V^*\cong V$.
A morphism $f : (V,\alpha)\to (W,\beta)$ in $\CF$ is an $\mathbb{F}$-linear function $f : V\otimes V\to W\otimes W$ for
which there exist $\lambda \in \mathbb{F}$ and $\mathbb{F}$-linear functions $a : V\to W$ and
$b : W \to V$ such that $f = a \otimes b^t$ where $ab=\lambda 1_W$, $ba=\lambda 1_V$,
and $b^t=\beta b^* \alpha^{-1}$. Composition is that of $\mathbb{F}$-linear functions. 
Each monoid $\CF(V,\alpha)$ is a monoid of Lie type; see Example 2.3 of \cite{OP1991}. 
By Remark~\ref{RmkKovacsId}, every $(V,\alpha)$ has a Kov\'acs idempotent. 
    \end{example}
     \begin{example}\label{KovacsSteinEx}
 In September 2022, Benjamin Steinberg told me that Itamar Stein had proved the semisimplicity
 of the $R$-algebra $R\Delta_{\bot}(\mathbf{a})$ for $R$ any field; see Example~\ref{exlocpo} for our notation. The proof has now appeared in \cite{Stein2022}.
    \end{example}
 
    \begin{lemma}\label{centralnature} 
 Assume the pigeon-hole principle and that each object $C$ has a Kov\'acs idempotent $\ell_C$. Then $\ell_C$ is a central idempotent in the $R$-algebra $R\CF(C)$
 and the morphisms $\ell_C : C \to C$ are the components of a natural endomorphism 
 of the identity functor of $R\CG$. 
 \end{lemma}
 \begin{proof}
 Take any $f : C\to C$ in $\CF$. Then, by the two-sided ideal property, both $\ell_Cf$ and $f\ell_C$ are in $R\CG'(C,C)$ so,
 by the identity element property, 
 $\ell_Cf = (\ell_Cf)\ell_C = \ell_C(f\ell_C) = f\ell_C$, proving the first clause.
 
 Now take $g \in \CG(C,D)$. For any $s\in \CG'(D,D)$, we have $g^{-1}sg\in\CG'(C,C)$. So $(g^{-1}sg)\ell_C = g^{-1}sg$.
 This shows that $s(g\ell_Cg^{-1}) = s$ yielding that $g\ell_Cg^{-1}$ is a right identity in $R\CG'(D,D)$. By uniqueness,
 $g\ell_Cg^{-1} = \ell_D$ which proves the claimed naturality.    
 \end{proof}    
   
\begin{proposition}
Assume that $\CF$ is stiff and that each object $C$ has a Kov\'acs idempotent $\ell_C$.
Then the idempotents $p_A=1_A -\ell_A : A\to A$ satisfy the Idempotent Axiom of Section~\ref{MainThmSection}.
    \end{proposition}
    \begin{proof}
 By duality, Lemma~\ref{centralnature} and $p'_A=\ell_A \in R\CG'(A,A)$, it remains to prove condition p1.
 Take $f\in \CF(A,B)$ with $f\notin \CM$. Then $f$ is not stiff so there is a factorization $f=vu$ with $u\in \CG'(A,A)$. Consequently, $fp_A = f(1_A -\ell_A)= f - vu\ell_A = f - vu = 0$.
 \end{proof} 

\section{Monoidal structure on $[\CF,\CV]$}

The bicategory $\mathfrak{R} = \CV\text{-}\mathrm{Mod}$ is autonomous symmetric monoidal (see \cite{60}).
The tensor product is the usual tensor product of $\CV$-enriched categories as per \cite{EilKel1966}:
take cartesian product of object sets and tensor over $R$ of hom $R$-modules.

Comonoidales (= pseudocomonoids) in $\mathfrak{R}$ are promonoidal $\CV$-categories in the sense of Day \cite{Day3}.
As a consequence of Theorem~\ref{mainthm}, any promonoidal structure on $R\CG$ transports across the adjoint equivalence $E\dashv M$ to a promonoidal structure on $R\CF$. In the case $\CF = \mathrm{vect}_{\mathbb{F}}$
for a finite field $\mathbb{F}$ and for $R=\mathbb{C}$ the complex numbers, a braided promonoidal structure
was defined on $R\CG$ in Remark 4.3 of \cite{49}. 
Here we will recall the general setting of Example 4 in Section 2 of \cite{113} and transport the promonoidal structure
to $R\CF$.  

{\em The assumptions of this section are that $\CF$ is an abelian category, $\CE$ consists of the epimorphisms, $\CM$ consists of the monomorphisms, and the Idempotent Axiom of Section~\ref{MainThmSection} holds.}

This promonoidal structure (unlike the braiding) actually exists on the ordinary category $\CG$ not just the $R$-linear category $R\CG$. It uses the
exact sequences of $\CF$. We have
\begin{align*}
P_{\CG} : \CG^{\mathrm{op}} \times \CG^{\mathrm{op}} \times \CG &\lra \mathrm{Set} 
\end{align*}
defined by 
$$
P_{\CG}(A,B;C) = \bigl\{(x,y) \bigm |  0 \to A \stackrel{x}{\to} C \stackrel{y}{\to} B \to 0
\text{ is a short exact sequence in } \CF  \bigr\} \ ,
$$
$P_{\CG}(g_1,g_2;g_3)(x,y) = (g_3xg_1,g_2^{-1}yg_3^{-1})$ for $g_1\in\CG(A',A), g_2\in\CG(B',B), g_3\in\CG(C,C')$,  and
\begin{align*}
J_{\CG} : \CG &\lra \mathrm{Set} 
\end{align*}
defined by
$$
J_{\CG}A = \begin{cases}
               1 & \text{for } A = 0 \\
\emptyset &\text{otherwise}.
\end{cases}
$$
The associativity isomorphism 
$$\alpha : \int^X P_{\CG} (X,C;D) \times P_{\CG}(A,B;X) \cong \int^Y P_{\CG}(A,Y;D) \times P_{\CG}(B,C;Y)$$ 
is defined by
$$
\alpha [(u,v),(x,y)] = [(x',y'),(u',v')] 
$$
as in the below $3 \times 3$ diagram whose rows and columns are short exact sequences in $\CF$. 
\begin{eqnarray*}
\xymatrix{
A \ar[d]_-{x} \ar[r]^-{1} & A \ar[d]_-{x'} \ar[r] & 0 \ar[d]  \\
 X \ar[d]_-{y} \ar[r]^-{u} & D \ar[d]_-{y'} \ar[r]^-{v} & C \ar[d]^-{1}  \\
 B \ar[r]_-{u'} & Y  \ar[r]_-{v'} & C   
}
\end{eqnarray*}

Day convolution yields the monoidal structure on $[\CG,\CV]$ with tensor product
\begin{eqnarray*}
(S\otimes_{\CG} T)C = \sum_{H\le C} {SH\ox T(C/H)}
\end{eqnarray*}
 where the sum runs over subobjects $H$ of $C$ in $\CF$ and $C/H$ is the quotient object.
 
 The promonoidal structure $RP_{\CG}$ on $R\CG$ transports across the adjoint equivalence of Theorem~\ref{mainthm}
 to a promonoidal structure $RP_{R\CF}$ on $R\CF$ defined by
 \begin{align}\label{promononRF}
P_{R\CF}(A,B;C) & = \int^{A',B',H}{M(H,C)\ox RP_{\CG}(A',B';H)\ox E(A,A')\ox E(B,B')} \nonumber \\
& \cong \int^{H}{M(H,C)\ox (E(A,-)\ox_{\CG} E(B,-))H} \nonumber \\
& \cong \sum_{K\le H\le C}{R\CE(A,jK)\ox R\CE(B,j(H/K))} \ .
\end{align}
 
 Day convolution yields the monoidal structure on $[\CF,\CV]$ with tensor product
\begin{eqnarray*}
(F\otimes_{\CG} G)C = \sum_{K\le H\le C} \int^{A,B\in \CF}{R\CE(A,jK)\ox R\CE(B,j(H/K))\ox FA\ox GB} \ .
\end{eqnarray*}
\begin{remark}
While the promonoidal structure \eqref{promononRF} on $R\CF$ is not the $R$-linearization of one on $\CF$,
it does exist at the level of the free category $p\CF$ enriched in pointed sets on $\CF$. This is clear from
the formula for the effect of $E$ on morphisms given in Lemma~\ref{tilderhosplitting}.  
\end{remark}

\end{document}